\documentclass[a4paper,11pt,twoside]{amsart}

\usepackage[latin1]{inputenc}
\usepackage[T1]{fontenc}
\usepackage{amsmath}
\usepackage{amsthm}
\usepackage{amssymb}
\usepackage[all]{xypic}

\newtheorem{thm}{Theorem}[section]
\newtheorem{prop}[thm]{Proposition}
\newtheorem{lem}[thm]{Lemma}
\newtheorem{cor}[thm]{Corollary}

\numberwithin{equation}{section}

\theoremstyle{definition}
\newtheorem{definition}[thm]{Definition}
\newtheorem{remark}[thm]{Remark}

\newcommand{\Db}{{\rm D}^{\rm b}}
\newcommand{\Km}{{\rm Km}}
\newcommand{\Eq}{{\rm Eq}}

\newcommand{\Br}{{\rm Br}}
\newcommand{\NS}{{\rm NS}}
\newcommand{\Pic}{{\rm Pic}}
\newcommand{\de}{{\rm d}}
\newcommand{\KS}{{\rm KS}}

\newcommand{\cat}[1]{\begin{bf}#1\end{bf}}
\newcommand{\coh}{\cat{Coh}}

\newcommand{\mor}[1][]{\xrightarrow{#1}}
\newcommand{\isomor}{\mor[\sim]}

\newcommand{\K}{{\rm K}}
\newcommand{\U}{{\rm U}}
\newcommand{\KE}{{\rm Ker}}

\newcommand{\cal}{\mathcal}

\newcommand{\ke}{{\cal E}}
\newcommand{\kf}{{\cal F}}
\newcommand{\kg}{{\cal G}}
\newcommand{\kh}{{\cal H}}

\newcommand{\ko}{{\cal O}}
\newcommand{\kp}{{\cal P}}

\newcommand{\ZZ}{\mathbb{Z}}
\newcommand{\QQ}{\mathbb{Q}}

\newcommand{\CC}{\mathbb{C}}

\begin{document}

\title[Derived categories and Kummer varieties]{Derived categories and Kummer varieties}

\author{Paolo Stellari}

\address{Dipartimento di Matematica ``F.
Enriques'', Universit{\`a} degli Studi di Milano, Via Cesare Saldini
50, 20133 Milano, Italy}

\email{stellari@mat.unimi.it}

\begin{abstract} We prove that if two abelian varieties have equivalent derived
categories then the derived categories of the smooth stacks
associated to the corresponding Kummer varieties are equivalent as
well. The second main result establishes necessary and sufficient
conditions for the existence of equivalences between the twisted
derived categories of two Kummer surfaces in terms of Hodge
isometries between the generalized transcendental lattices of the
corresponding abelian surfaces.\end{abstract}

\maketitle

\section{Introduction}

The {\it Kummer variety} of an abelian variety $A$ is the quotient
$\K(A):=A/\langle\iota\rangle$, where $\iota(a)=-a$ for any $a\in
A$. The singular variety $\K(A)$ has an orbifold structure and
it admits a minimal crepant resolution $\Km(A)$
if and only if the dimension of $A$ is 2. In this case $\Km(A)$ is a
K3 surface (i.e.\ it is simply connected and its dualizing sheaf is
trivial) and it is called the {\it Kummer surface of $A$}. More generally, we can associate to the global quotient $\K(A)$ the smooth quotient stack
$[A/\langle\iota\rangle]$.

In \cite{HLOY} Hosono, Lian, Oguiso and Yau proved that,
\begin{itemize}
\item[(A)] {\it given two abelian surfaces $A$ and $B$, $\Db(A)\cong\Db(B)$ if
and only if $$\Db(\Km(A))\cong\Db(\Km(B)).$$}
\end{itemize}
Their argument runs as follows: They notice that, due to the
geometric construction of the Kummer surfaces $\Km(A)$ and $\Km(B)$,
the transcendental lattices of $A$ and $B$ are Hodge-isometric if
and only if the transcendental lattices of $\Km(A)$ and $\Km(B)$ are
Hodge-isometric. Then, they apply a deep result of Orlov which says
that two abelian or K3 surfaces have equivalent derived categories
if and only if their transcendental lattices are Hodge-isometric
(see Theorem \ref{thm:orlov1}). From this it is evident that (A) can
be reformulated in the following way:
\begin{itemize}
\item[(B)] {\it given two abelian surfaces $A$ and $B$,
$\Db(\Km(A))\cong\Db(\Km(B))$ if and only if there exists a Hodge
isometry between the transcendental lattices of $A$ and $B$.}
\end{itemize}
Since Mukai proved in \cite{Mu} that two K3 surfaces with Picard
number greater than 11 and with Hodge-isometric transcendental
lattices are isomorphic, (A) and (B) are equivalent to the
following statement:
\begin{itemize}
\item[(C)] {\it given two abelian surfaces $A$ and $B$,
$\Db(A)\cong\Db(B)$ if and only if $\Km(A)\cong\Km(B).$}
\end{itemize}

The aim of this paper is to address (A), (B) and (C) in two more
general contexts. Our first result shows that if $A_1$ and $A_2$
are abelian varieties with equivalent derived categories, then the
derived categories of the stacks $[A_1/\langle\iota\rangle]$ and
$[A_2/\langle\iota\rangle]$ are equivalent as well. In fact we
will prove the following:

\begin{thm}\label{cor:equiv} Let $A_1$ and $A_2$ be abelian varieties.
If $\Db(A_1)\cong\Db(A_2)$, then there exists a Fourier-Mukai
equivalence
$\Db([A_1/\langle\iota\rangle])\cong\Db([A_2/\langle\iota\rangle])$.

Conversely, if $\Db([A_1/\langle\iota\rangle])$ and
$\Db([A_2/\langle\iota\rangle])$ are equivalent, then there is an
isomorphism of Hodge structures $\widetilde
H(A_1,\QQ)\cong\widetilde H(A_2,\QQ)$.\end{thm}

The Hodge structures mentioned in the second part of the previous
statement will be defined in Section \ref{subsec:Hodge} and the
proof of this result will occupy almost all Section
\ref{sec:equivav}. As we will show in Section
\ref{subsec:geomappl1}, when we deal with abelian surfaces this
result leads to a direct proof of one implication in (A). An
application to the number of birational generalized Kummer varieties
is given in Proposition \ref{cor:genKum}.

Our second main result treats the two-dimensional twisted case.
Indeed, according to (B), we prove that the twisted derived
categories of two Kummer surfaces are equivalent if and only if the
generalized transcendental lattices of the corresponding abelian
surfaces are Hodge isometric. More precisely the result (proved in
Sections \ref{subsec:Br} and \ref{subsec:proofsec}) is as follows:

\begin{thm}\label{thm:main2} Let $A_1$ and $A_2$ be abelian
surfaces. Then the following two conditions are equivalent:
\begin{itemize}
\item[(i)] there exist $\alpha_i$ in the Brauer group of $\Km(A_i)$
and an equivalence between the derived categories
$\Db(\Km(A_1),\alpha_1)$ and $\Db(\Km(A_2),\alpha_2)$;
\item[(ii)] there exist $\beta_i$ in the Brauer group of $A_i$ such
that the twisted abelian surfaces $(A_1,\beta_1)$ and
$(A_2,\beta_2)$ have Hodge-isometric generalized transcendental
lattices.
\end{itemize}
Furthermore, if one of these two equivalent conditions holds true,
then $A_1$ and $A_2$ are isogenous.\end{thm}

The notations an definitions involved in the formulation of the
previous result will be explained in Sections \ref{subsec:Hodge} and
\ref{subsec:Br}. We will observe that the analogues of (A) and (C)
in the twisted setting are no longer true (see Remark
\ref{rmk:commenti}). Nevertheless we completely generalize the
results in \cite{HLOY} about the number of Kummer structures on K3
surfaces in the twisted context (Proposition \ref{prop:twnum}). A
geometric example involving abelian surfaces with Picard number two
is discussed.

\section{Derived categories of abelian varieties and K3 surfaces}\label{sec:prelim}

In this section we recall some facts and definitions concerning the
derived categories of coherent sheaves on abelian varieties and K3
surfaces. In the following pages $\Db(X)$ will always mean the
bounded derived category of coherent sheaves on the smooth
projective variety $X$ (we will also use the same notation for the
bounded derived category of coherent sheaves on a smooth quotient
stack according to \cite{Ka}).

Suppose that $X_1$ and $X_2$ are smooth projective varieties. Let
$\Db(X_1)$ and $\Db(X_2)$ be the bounded derived categories of
coherent sheaves on $X_1$ and $X_2$. Orlov proved in \cite{Or1} that
any equivalence $\Phi:\Db(X_1)\rightarrow\Db(X_2)$ is a {\it
Fourier-Mukai equivalence}, i.e. there exists $\ke\in\Db(X_1\times
X_2)$ and an isomorphism of functors
\begin{eqnarray}\label{eqn:FMtr}
\Phi\cong{\bf R}p_{2*}(\ke\stackrel{\bf L}{\otimes}p_1^*(-)),
\end{eqnarray}
where $p_i:X_1\times X_2\rightarrow X_i$ is the projection and
$i\in\{1,2\}$ (see also \cite{CS} for a more general statement). The
complex $\ke$ is the {\it kernel} of $\Phi$ and it is uniquely (up
to isomorphism) determined. We write $\Phi_\ke$ for a Fourier-Mukai
equivalence whose kernel is $\ke$. In general, given
$\ke\in\Db(X_1\times X_2)$, we write $\Phi_\ke$ for a functor
defined as in \eqref{eqn:FMtr} (notice that $\Phi_\ke$ is not
necessarily an equivalence).

\subsection{Derived categories of abelian varieties} Assume that $A_1$ and $A_2$ are abelian varieties of dimension
$d$. For $i\in\{1,2\}$, let $\kp_i$ be the Poicar{\'e} line bundle
on $A_i\times\widehat{A_i}$, let $\mu_i:A_i\times A_i\rightarrow
A_i\times A_i$ be the isomorphism such that $(a,b)\mapsto (a+b,b)$
and let $\Phi_i:=\mu_{i*}\circ(\mathrm{id}\times\Phi_{\kp_i})$. If
$\Phi_\ke:\Db(A_1)\rightarrow\Db(A_2)$ is a Fourier-Mukai
equivalence with kernel $\ke$, we get the following commutative
diagram:
\begin{eqnarray}\label{eqn:der}
\xymatrix{\Db(A_1\times\widehat{A_1})\ar[dr]^{\mathrm{id}\times\Phi_{\kp_1}}\ar[rrr]^{F_\ke}\ar[dd]_{\Phi_1}&&&
\Db(A_2\times\widehat{A_2})\ar[dl]_{\mathrm{id}\times\Phi_{\kp_2}}\ar[dd]_{\Phi_2}\\
&\Db(A_1\times A_1)\ar[dl]_{\mu_{1*}}& \Db(A_2\times A_2)\ar[dr]^{\mu_{2*}}&\\
\Db(A_1\times
A_1)\ar[rrr]^{\Phi_\ke\times\Phi_{\ke_R}}&&&\Db(A_2\times A_2),}
\end{eqnarray}
where $F_\ke$ is the functor completing the diagram,
$\ke_R=\ke^\vee[d]$ and $\Phi_\ke\times\Phi_{\ke_R}$ is the
Fourier-Mukai equivalence whose kernel is $\ke\boxtimes\ke_R$.
Observe that since $\Phi_\ke$, $\Phi_1$ and $\Phi_2$ are
equivalences, $\Phi_\ke\times\Phi_{\ke_R}$ and $F_\ke$ are
equivalences as well.

For $i\in\{1,2\}$, the K\"{u}nneth formula yields a decomposition
$$H_1(A_i\times\widehat{A_i},\ZZ)\cong H_1(A_i,\ZZ)\oplus
H_1(\widehat{A_i},\ZZ).$$ Since $H_1(\widehat{A_i},\ZZ)\cong
H_1(A_i,\ZZ)^\vee$, the group $H_1(A_i\times\widehat{A_i},\ZZ)$ is
endowed with a natural quadratic form. Indeed, if
$(a_1,\alpha_1),(a_2,\alpha_2)\in H_1(A_i\times\widehat{A_i},\ZZ)$,
we define
\[
\langle(a_1,\alpha_1),(a_2,\alpha_2)\rangle_i:=\alpha_1(a_2)+\alpha_2(a_1),
\]
where $i\in\{1,2\}$. Consider the set of isomorphisms
\[
\U(A_1,A_2):=\{f\in\mathrm{Isom}(A_1\times\widehat{A_1},A_2\times\widehat{A_2}):\langle
f_*(a_1,\alpha_1),f_*(a_2,\alpha_2)\rangle_2=\langle
(a_1,\alpha_1),(a_2,\alpha_2)\rangle_1\}.
\]

\begin{thm}\label{thm:orlov} {\bf(\cite{Or}, Theorem 2.19 and Proposition 4.12.)}  Let $A_1$ and $A_2$ be abelian varieties. If $\Phi_\ke:\Db(A_1)\rightarrow\Db(A_2)$ is an
equivalence, then, for any $\kf\in\Db(A_1)$,
$$F_\ke(\kf)=f_{\ke\,*}(\kf)\otimes N_\ke,$$ where $F_\ke$ is the
equivalence in \eqref{eqn:der}, $f_\ke\in\U(A_1,A_2)$ and
$N_\ke\in\Pic(A_2\times\widehat{A_2})$. Moreover, there exists a
surjective map
\begin{eqnarray}\label{eqn:gamma}
\gamma:\Eq(\Db(A_1),\Db(A_2))\longrightarrow\U(A_1,A_2)
\end{eqnarray}
such that $\gamma(\Phi_\ke)=f_\ke$, where $\Eq(\Db(A_1),\Db(A_2))$
is the set of equivalences between $\Db(A_1)$ and
$\Db(A_2)$.\end{thm}

\subsection{Hodge structures and derived categories}\label{subsec:Hodge} If $X$ is a smooth
projective variety of dimension $d$, we denote by $\widetilde
H(X,\QQ)$ the even cohomology group $H^{2*}(X,\QQ)$ with the
weight-$d$ Hodge structure defined as follows:
\begin{eqnarray}\label{eqn:HSnuova}
\widetilde{H}^{p,q}(X)=\bigoplus_{p-q=r-s} H^{r,s}(X),
\end{eqnarray}
where $H^{r,s}(X)$ is the $(r,s)$-part of the usual Hodge
decomposition of $H^{r+s}(X,\CC)\subset H^{2*}(X,\CC)$. An
equivalent way to put on $H^*(X,\CC)$ such a Hodge structure could
be obtained considering the natural grading on the Hochschild
homology of $X$ (see, for example, \cite{C3}).

Suppose now that $X$ is either an abelian or a K3 surface,
$H^{2,0}(X)=\langle\sigma_X\rangle$ and $B$ is any class in
$H^2(X,\QQ)$. Then
$$\varphi:=\exp(B)(\sigma_X)=\sigma_X+B\wedge\sigma_X\in
H^2(X,\CC)\oplus H^4(X,\CC)$$ is a {\it generalized Calabi-Yau
structure on $X$} (for a complete picture see \cite{Hu}). Let
$T(X,B)$ be the minimal primitive sublattice of $H^2(X,\ZZ)\oplus
H^4(X,\ZZ)$ such that $\varphi\in T(X,B)\otimes\CC$. The lattice
$T(X,B)$ is the {\it generalized transcendental lattice of
$\varphi$} (see \cite{Hu} and \cite{HS1}). Let $\widetilde H(X,\ZZ)$
be the $\ZZ$-module $H^0(X,\ZZ)\oplus H^2(X,\ZZ)\oplus H^4(X,\ZZ)$
endowed with the \emph{Mukai pairing}
\[
\langle(a_0,a_2,a_4),(b_0,b_2,b_4)\rangle=a_2\cdot b_2-a_0\cdot
b_4-a_4\cdot b_0,
\]
where $(a_0,a_2,a_4),(b_0,b_2,b_4)\in H^{2*}(X,\ZZ)$ and
``$\cdot$'' is the cup-product. We write $\widetilde H(X,B,\ZZ)$
for the lattice $\widetilde H(X,\ZZ)$ with the weight-two Hodge
structure such that
\[
\widetilde H^{2,0}(X,B):=\exp(B)(\widetilde H^{2,0}(X))
\]
and $\widetilde H^{1,1}(X,B)$ is its orthogonal complement in
$H^2(X,\CC)$. It is clear that $T(X,B)$ inherits from $\widetilde
H(X,B,\ZZ)$ a weight-two Hodge structure. By definition,
$T(X)=T(X,0)$ is the {\it transcendental lattice of $X$} and
$$\NS(X):=T(X)^\perp\subset H^2(X,\ZZ)$$ is the {\it
N{\'e}ron-Severi group of $X$}. The number
$\rho(X):=\mathrm{rk}\NS(X)$ is the {\it Picard number} of $X$. If
$L_1$ and $L_2$ are lattices endowed with a weight-$k$ Hodge
structure, then an isometry $f:L_1\rightarrow L_2$ is a {\it Hodge
isometry} if it preserves the Hodge structures.

For abelian and K3 surfaces, Orlov proved in \cite{Or1} (using
results in \cite{Mu}) the following theorem:

\begin{thm}\label{thm:orlov1} {\bf(\cite{Or1}, Theorem 3.3.)} Let $X_1$ and $X_2$ be either
abelian or K3 surfaces. Then the following two conditions are
equivalent:

{\rm (i)} there exists an equivalence $\Db(X_1)\cong\Db(X_2)$;

{\rm (ii)} there exists a Hodge isometry $T(X_1)\cong
T(X_2)$.\end{thm}

\section{Derived categories of the smooth stacks}\label{sec:equivav}

This section is mainly devoted to the proof of Theorem
\ref{cor:equiv}. As it will turn out, such a proof, which will be
given in Section \ref{subsec:proof}, relies on some results about
the equivariant derived categories of coherent sheaves on abelian
varieties proved in Section \ref{subsec:equivcat}. In Section
\ref{subsec:geomappl1} some geometric applications are discussed.

\subsection{Equivariant derived categories and abelian varieties}\label{subsec:equivcat}
Consider the simple case of an abelian variety $A$ with the action
of $G:=\ZZ/2\ZZ$ induced by the automorphism $\iota:A\rightarrow A$
such that $\iota(a)=-a$, for any $a\in A$. A $G$-linearization for a
coherent sheaf $\ke\in\coh(A)$ is an isomorphism
$\lambda:\ke\rightarrow\iota^*\ke$ such that
$\iota^*(\lambda)=\lambda$ and
$\iota^*(\lambda)\circ\lambda=\lambda\circ\lambda=\mathrm{id}$.

$\coh^G(A)$ is the abelian category whose objects are the pairs
$(\ke,\lambda)$, where $\ke\in\coh(A)$ admits a $G$-linearization
and $\lambda$ is a $G$-linearization for $\ke$. The morphisms in
$\coh^G(A)$ are the morphisms in $\coh(A)$ compatible with the
$G$-linearizations. We define $\Db_G(A):=\Db(\coh^G(A))$ to be the
bounded derived category of $\coh^G(A)$.  A complete discussion
about the general case when $G$ is any finite group acting on a
smooth projective variety can be found in \cite{BKR}.

If $A_1$ and $A_2$ are abelian varieties and $G_\Delta\cong\ZZ/2\ZZ$
acts on $A_1\times A_2$ via the automorphism $\iota\times\iota$, the
set of $G_\Delta$-invariant equivalences has the following
description:
\[
\Eq(\Db(A_1),\Db(A_2))^{G_{\Delta}}=\{\Phi_\kg\in\Eq(\Db(A_1),\Db(A_2)):\kg\in\Db(A_1\times
A_2)\mbox{ is $G_\Delta$-invariant}\}.\\
\]
An equivalence $\Phi:\Db_G(A_1)\cong\Db_G(A_2)$ is a Fourier-Mukai
equivalence if there is an isomorphism as in \eqref{eqn:FMtr}, where
the kernel $\ke$ is in $\Db_{G\times G}(A_1\times A_2)$.
$\Eq(\Db_G(A_1),\Db_G(A_2))$ is the set whose elements are the
equivalences of this type.

\begin{prop}\label{thm:main1} Let $A_1$ and $A_2$ be abelian varieties
and let $G=\ZZ/2\ZZ$ act on $A_1$ and $A_2$ as above. Then the
restriction
\[
\gamma:\Eq(\Db(A_1),\Db(A_2))^{G_\Delta}\longrightarrow\U(A_1,A_2)
\]
of the map in \eqref{eqn:gamma} is surjective and
$\Eq(\Db_G(A_1),\Db_G(A_2))$ is non-empty if $\U(A_1,A_2)$ is
non-empty.\end{prop}

\begin{proof}  By definition, we can think of any $f\in\U(A_1,A_2)$ as represented by a matrix
\[
\left(\begin{array}{cc} x_f & y_f\\ z_f & w_f\end{array}\right).
\]
Define $S(A_1,A_2):=\{f\in\U(A_1,A_2): y_f\mbox{ is an isogeny}\}$
and let $f\in S(A_1,A_2)$. Using results from \cite{Mu1}, Orlov
proved in \cite{Or} (see, in particular, \cite[Prop.\ 4.12]{Or})
that there exists a vector bundle $\ke$ on $A_1\times A_2$ with the
following properties:
\begin{itemize}
\item[(a)] $\ke$ is simple and $\Phi_\ke$ is an equivalence;
\item[(b)] for any $(a,b)\in A_1\times A_2$, if $T_{(a,b)}$ is the translation with respect to the point $(a,b)$,
then $T_{(a,b)*}\ke\cong\ke\otimes P$ for some
$P\in\mathrm{Pic}^0(A_1\times A_2)$;
\item[(c)] $\gamma(\Phi_\ke)=f$.
\end{itemize}
Consider now the sheaf $\kf:=(\iota\times\iota)^*\ke$. It is clear
that $\gamma(\Phi_\kf)=\gamma(\Phi_\ke)=f$.

For a brief proof of this fact, consider the maps $\Phi_i$,
$\Phi_{\kp_i}$ and $\mu_i$ in \eqref{eqn:der}. A straightforward
calculation shows that $(\iota\times\iota)^*\kp_i\cong\kp_i$.
Moreover, $\mu_i$ is a morphism of abelian varieties. Hence
$(\iota\times\iota)^*\circ(\mathrm{id}\times\Phi_{\kp_i})\circ(\iota\times\iota)^*=\mathrm{id}\times\Phi_{\kp_i}$
and
$(\iota\times\iota)^*\circ\mu_{i*}\circ(\iota\times\iota)^*=\mu_{i*}$.
This implies that
$(\iota\times\iota)^*\circ\Phi_i\circ(\iota\times\iota)^*=\Phi_i$,
for $i\in\{1,2\}$. Since $\Phi_\kf=\iota^*\circ\Phi_\ke\circ\iota^*$
and
$\Phi_\kf\times\Phi_{\kf_R}=(\iota\times\iota)^*\circ(\Phi_\ke\times\Phi_{\ke_R})\circ(\iota\times\iota)^*$,
we rewrite the commutative diagram \eqref{eqn:der} in the following
way:
\begin{eqnarray}\label{eqn:dia1}
\xymatrix{\Db(A_1\times\widehat{A_1})\ar[d]_{\Phi_1}\ar[r]_{(\iota\times\iota)^*}
\ar@/^2pc/[rrr]^{F_\kf} &
\Db(A_1\times\widehat{A_1})\ar[d]_{\Phi_1}\ar[r]^{F_\ke} &
\Db(A_2\times
\widehat{A_2})\ar[d]_{\Phi_2}\ar[r]_{(\iota\times\iota)^*} &
\Db(A_2\times\widehat{A_2})\ar[d]_{\Phi_2}\\
\Db(A_1\times A_1)\ar[r]^{(\iota\times\iota)^*}
\ar@/_2pc/[rrr]^{\Phi_\kf\times\Phi_{\kf_R}}
 & \Db(A_1\times A_1)\ar[r]^{\Phi_\ke\times\Phi_{\ke_R}} & \Db(A_2\times
 A_2)\ar[r]^{(\iota\times\iota)^*} & \Db(A_2\times A_2).}
\end{eqnarray}
By Theorem \ref{thm:orlov}, for any
$\kg\in\Db(A_1\times\widehat{A_1})$,
$F_\kf(\kg)=f_{\kf*}(\kg)\otimes N_\kf$, for some
$f_\kf\in\U(A_1,A_2)$ and $N_\kf\in\Pic(A_2\times\widehat{A_2})$.
Hence, from \eqref{eqn:dia1} we deduce that
{\setlength\arraycolsep{2pt}
\begin{eqnarray*}
F_\kf(\kg)&=&((\iota\times\iota)^*\circ
F_\ke\circ(\iota\times\iota)^*)(\kg)\\&=&((\iota\times\iota)^*\circ
f_{\ke*}\circ(\iota\times\iota)^*)(\kg)\otimes
M\\&=&f_{\ke*}(\kg)\otimes M,
\end{eqnarray*}}for some $M\in\Pic(A_2\times\widehat{A_2})$. Observe that the last
equality holds true because $f_\ke$ is a morphism of abelian
varieties. This proves that
$\gamma(\Phi_\kf)=\gamma(\iota^*\circ\Phi_\ke\circ\iota^*)=\gamma(\Phi_\ke)$
which is what we claimed.

Due to this last remark and to \cite[Cor.\ 3.4]{Or}, there exist
$a\in A_1$ and $\alpha\in\widehat{A_1}$ such that
\begin{eqnarray}\label{eqn:orlsh}
\kf=T_{(a,0)*}\ke\otimes p^*P_\alpha[i],
\end{eqnarray}
where $p:A_1\times A_2\rightarrow A_1$ is the projection, $i$ is an
integer and $P_\alpha$ is the degree zero line bundle on $A_1$
corresponding to $\alpha$. In the following arguments, without loss
of generality, we will forget about the shift $[i]$ in
\eqref{eqn:orlsh}.

Since $\ke$ satisfies (b), from \eqref{eqn:orlsh} we get
$\kf=\ke\otimes Q$, where $Q$ is a degree zero line bundle on
$A_1\times A_2$. Let $N\in\mathrm{Pic}^0(A_1\times A_2)$ be such
that $N^2=Q$ and consider the sheaf $\ke_f:=\ke\otimes N$. It is
easy to see that
{\setlength\arraycolsep{2pt}
\begin{eqnarray*}
(\iota\times\iota)^*(\ke_f) & =&(\iota\times\iota)^*(\ke\otimes
N)\nonumber \\
 & \cong&\ke\otimes Q\otimes N^\vee\\
 & \cong&\ke\otimes N\\
 & =&\ke_f.\nonumber
\end{eqnarray*}}Due to \cite[Prop.\ 3.3]{Or} and to (c),
$\gamma(\Phi_{\ke_f})=\gamma(\Phi_\ke)=f$.

Let $f\in\U(A_1,A_2)$. Orlov observed in \cite[Sect.\ 4]{Or} that
there exist $g_1\in S(A_1,A_2)$ and $g_2\in S(A_2,A_2)$ such that
$f=g_2\circ g_1$. From its very definition, the map $\gamma$ in
Theorem \ref{thm:orlov} preserves the compositions. Hence $\gamma$
restricts to a surjective map
$\gamma:\Eq(\Db(A_1),\Db(A_2))^{G_\Delta}\rightarrow\U(A_1,A_2)$.

To prove the second claim in Proposition \ref{thm:main1}, consider
the set
\[
\KE(A_1,A_2,G_\Delta):=\{(\kg,\lambda)\in\Db_{G_\Delta}(A_1\times
A_2):\Phi_\kg\in\Eq(\Db(A_1),\Db(A_2))\}.
\]
Since the group cohomology
$H^2(\mathbb{Z}/2\mathbb{Z},\mathbb{C}^*)$ is trivial, \cite[Thm.\
6]{Pl} shows the existence of two maps
\[
\begin{array}{l}
\psi_1:\KE(A_1,A_2,G_\Delta)\longrightarrow\Eq(\Db(A_1),\Db(A_2))^{G_{\Delta}}\\
\psi_2:\KE(A_1,A_2,G_\Delta)\longrightarrow\Eq(\Db_G(A_1),\Db_G(A_2))
\end{array}
\]
such that, for any $(\kg,\lambda)\in\KE(A_1,A_2,G_\Delta)$,
$\psi_1((\kg,\lambda))=\Phi_\kg$ and
$\psi_2((\kg,\lambda))=\Phi_\kh$, where
$\kh:=(\kg\oplus(\iota,\mathrm{id})^*\kg,\lambda')$ and $\lambda'$
is the natural $(G\times G)$-linearization induced by $\lambda$.

We previously proved that for any $f\in S(A_1,A_2)$, there exists
$\Phi_{\ke_f}\in\Eq(\Db(A_1),\Db(A_2))^{G_{\Delta}}$ such that
$\gamma(\Phi_{\ke_f})=f$. From \cite{Pl} it follows that $\psi_1$
is surjective and that the set $\KE(A_1,A_2,G_\Delta)$ is
non-empty if $\Eq(\Db(A_1),\Db(A_2))^{G_{\Delta}}$ is non-empty.
Hence, there exists $\Psi_f\in\KE(A_1,A_2,G_\Delta)$ such that
$\psi_1(\Psi_f)=\Phi_{\ke_f}$. The functor $\psi_2(\Psi_f)$ is in
$\Eq(\Db_G(A_1),\Db_G(A_2))$.\end{proof}

The special case $A_1=A_2$ is also treated in \cite{Pl}.

\subsection{Proof of Theorem \ref{cor:equiv}}\label{subsec:proof} Let $A_1$ and $A_2$
be abelian varieties and suppose that $\Db(A_1)\cong\Db(A_2)$. Due
to Theorem \ref{thm:orlov}, the set $\U(A_1,A_2)$ is non-empty.
Therefore, if $G=\ZZ/2\ZZ$ acts on $A_1$ and $A_2$ as prescribed
at the beginning of Section \ref{subsec:equivcat}, then
Proposition \ref{thm:main1} yields an equivalence
$\Psi:\Db_{G}(A_1)\isomor\Db_{G}(A_2)$.

Consider the stacks $[A_1/G]$ and $[A_2/G]$ (see \cite{Go} and
\cite{Ka}). For any $i\in\{1,2\}$, let $\Db([A_i/G])$ be the bounded
derived category of the abelian category $\coh([A_i/G])$ of coherent
sheaves on $[A_i/G]$ (see \cite{Ka}). Obviously
$\Db([A_i/G])\cong\Db_G(A_i)$, because
$\coh([A_i/G])\cong\coh^G(A_i)$, for any $i\in\{1,2\}$. This implies
that $\Psi$ can be rewritten as
$\Phi:\Db([A_1/G])\isomor\Db([A_2/G])$. Due to \cite{Ka}, $\Phi$ is
of Fourier-Mukai type (i.e.\ there is an isomorphism as in
\eqref{eqn:FMtr}). Hence, the first part of Theorem \ref{cor:equiv}
is proved.

Assume that an equivalence $\Phi:\Db([A_1/G])\isomor\Db([A_2/G])$
is given. As before, the results in \cite{Ka} imply that we can
think of $\Phi$ as a Fourier-Mukai equivalence whose kernel is a
$(G\times G)$-linearized complex $(\ke,\lambda)$. Obviously, the
inverse $\Phi^{-1}$ is a Fourier-Mukai equivalence as well.
Suppose that its kernel is $(\kf,\lambda')$. It is an easy
exercise to show that the kernel of the identity functor
$\mathrm{id}=\Phi\circ\Phi^{-1}:\Db_G(A_i)\rightarrow\Db_G(A_i)$
is the $(G\times G)$-linearized sheaf $(\ko_\Delta\oplus
(\iota,\mathrm{id})^*\ko_\Delta,\mu)$, where $\mu$ is the natural
linearization and $\Delta\hookrightarrow A_i\times A_i$ is the
diagonal embedding.

Consider the functors $\Phi_\ke$, $\Phi_\kf$ and $\Phi_{\ko_\Delta\oplus
(\iota,\mathrm{id})^*\ko_\Delta}$. Although they are
no longer equivalences, they induce the commutative diagram
\begin{eqnarray}\label{eqn:import}
\xymatrix{
\Db(A_1)\ar[d]_{\mathrm{ch}}\ar[r]^{\Phi_\ke}\ar@/^2pc/[rr]^{\Phi_{\ko_\Delta\oplus
(\iota,\mathrm{id})^*\ko_\Delta}} &
\Db(A_2)\ar[d]_{\mathrm{ch}}\ar[r]^{\Phi_\kf} &
\Db(A_1)\ar[d]_{\mathrm{ch}}\\
H^{2*}(A_1,\QQ)\ar[r]^{\Phi_\ke^H}
\ar@/_2pc/[rr]_{\Phi_{\ko_\Delta\oplus
(\iota,\mathrm{id})^*\ko_\Delta}^H}
 & H^{2*}(A_2,\QQ)\ar[r]^{\Phi_\kf^H} & H^{2*}(A_1,\QQ),}
\end{eqnarray}
where $\Phi_\ke^H:H^{2*}(A_1,\QQ)\rightarrow H^{2*}(A_2,\QQ)$ is
such that $\Phi^H_\ke(a)=p_{2*}(\mathrm{ch}(\ke)\cdot p_1^*(a))$ and
$p_i:A_1\times A_2\rightarrow A_i$ is the projection. Take analogous
definitions for $\Phi_\kf^H$ and $\Phi_{\ko_\Delta\oplus
(\iota,\mathrm{id})^*\ko_\Delta}^H$.

Observe that $(\iota,\mathrm{id})^*\ko_\Delta$ is the kernel of the
Fourier-Mukai equivalence $\iota^*:\Db(A_i)\isomor\Db(A_i)$. Since
$\iota^*$ acts as the identity on the cohomology lattice $\widetilde
H(A_i,\ZZ)$, from \eqref{eqn:import} we deduce
$$\Phi_\kf^H\circ\Phi_\ke^H=(\Phi_{\ko_\Delta\oplus
(\iota,\mathrm{id})^*\ko_\Delta})^H=2\mathrm{id}.$$ Hence
$\Phi_\ke^H$ is injective. Exchanging the roles of $\Phi_\ke$ and
$\Phi_\kf$ in \eqref{eqn:import}, we see that $\Phi_\ke^H$ is an
isomorphism of $\QQ$-vector spaces. In particular,
$\mathrm{dim}(A_1)=\mathrm{dim}(A_2)=n$.

The fact that the Hodge structures defined in \eqref{eqn:HSnuova}
are preserved follows from the standard argument for Fourier-Mukai
equivalences (see \cite[Prop.\ 5.39]{Hu1}). Indeed, one just needs
to observe that $\mathrm{ch}(\ke)\in\widetilde H^{2n,2n}(A_1\times
A_2)$. This concludes the proof of Theorem \ref{cor:equiv}.

\begin{remark} Of course, in general, $\Phi_\ke^H$ does not preserve the Mukai
pairing naturally defined on $H^{2*}(A_i,\QQ)$ by means of the cup
product (\cite[Chapter 5]{Hu1}). Indeed, it is easy to see that the
Mukai pairing is preserved up to a factor 2.\end{remark}

\subsection{Geometric applications}\label{subsec:geomappl1} Assume that $A_1$ and $A_2$
are abelian surfaces. The main result in \cite{BKR} yields an
equivalence
$\Psi_i:\Db([A_i/\langle\iota\rangle])\isomor\Db(\Km(A_i))$, for any
$i\in\{1,2\}$. Thus, if $\Upsilon:\Db(A_1)\isomor\Db(A_2)$ is a
Fourier-Mukai equivalence, we immediately get a second Fourier-Mukai
equivalence
\begin{eqnarray}\label{eqn:BKR1}
\Phi:\Db(\Km(A_1))\stackrel{\Psi_1^{-1}}{\longrightarrow}\Db([A_1/\langle\iota\rangle])\isomor\Db([A_2/\langle\iota\rangle])\stackrel{\Psi_2}{\longrightarrow}\Db(\Km(A_2)),
\end{eqnarray}
where the middle equivalence is produced by Theorem \ref{cor:equiv}
and the kernel of $\Phi$ can be easily computed using \cite{BKR}.
This leads to a different and explicit proof of the ``only if''
implication in (A) without using the lattice theoretical description
of the transcendental lattices of an abelian surface and of the
associated Kummer surface.

Let us discuss a second geometric application. Assume that $A$ is
an abelian surface. We denote by $\mathrm{K}^n(A)$ the {\it $n$-th
generalized Kummer variety of $A$}. Recalling the construction in
\cite{Be}, we see that $\mathrm{K}^n(A)$ is the fiber over $0$
with respect to the map $\Psi$ which is the composition of the
morphisms in the following diagram:
\[
\Psi:\mathrm{Hilb}^{n+1}(A)\stackrel{\rho}{\longrightarrow}
\mathrm{Sym}^{n+1}(A)\stackrel{\sigma}{\longrightarrow} A,
\]
where $\rho$ is the Hilbert-Chow morphism and
$\sigma(a_1,\ldots,a_{n+1})=a_1+\ldots+a_{n+1}$. It is easy to
see that $\K^n(A)$ is smooth and that $\mathrm{K}^1(A)=\Km(A)$.
Furthermore, in \cite{Be} Beauville proved that these varieties are
examples of irreducible symplectic manifolds.

\begin{prop}\label{cor:genKum} Let $A$ be an abelian surface and let
$n\geq 2$ be an integer. The number of generalized Kummer
varieties $\K^n(B)$ birational to $\K^n(A)$ is finite up to
isomorphisms. Moreover if $\K^n(A)$ and $\K^n(B)$ are birational,
then $\Km(A)\cong\Km(B)$ and $A$ and $B$ are isogenous.\end{prop}

\begin{proof} Let $A_1$ and $A_2$ be abelian surfaces and let $\varphi$ be
a birational morphism between $\mathrm{K}^n(A_1)$ and
$\mathrm{K}^n(A_2)$. Obviously, $\varphi$ induces an isomorphism
$g:H^2(\mathrm{K}^n(A_1),\mathbb{Z})\isomor
H^2(\mathrm{K}^n(A_2),\mathbb{Z})$. Furthermore, there exists an
isometry of lattices $H^2(\mathrm{K}^n(A_i),\mathbb{Z})\isomor
H^2(A_i,\mathbb{Z})\oplus\mathbb{Z}[E_i]$, where $E_i$ is the
restriction to $\K^n(A)$ of the exceptional locus of
$\mathrm{Hilb}^{n+1}(A_i)$. The left hand side of the isomorphism is
endowed with the Beauville-Bogomolov form  while the quadratic form
on $H^2(A_i,\mathbb{Z})$ is the cup-product (see \cite[Lemma
4.10]{Be} and \cite[Prop.\ 4.11]{Y}).

Since $E_1$ and $E_2$ are algebraic, $g$ yields an isomorphism
$T(A_1)\cong T(A_2)$. Using Theorem \ref{thm:orlov1}, we get an
equivalence $\Db(A_1)\cong\Db(A_2)$. To prove that $A_1$ is
isogenous to $A_2$ observe that if $\Db(A_1)\cong\Db(A_2)$, then
$A_1\times\widehat{A_1}\cong A_2\times\widehat{A_2}$ (Theorem
\ref{thm:orlov}). Hence $A_1\times A_1$ and $A_2\times A_2$ are
isogenous and $A_1$ and $A_2$ are isogenous as well. On the other
hand, as there are only finitely many isomorphism classes of abelian
surfaces $A$ such that $\Db(A)\cong\Db(A_1)$ (see \cite[Prop.\
5.3]{BM}), the number of generalized Kummer varieties $\K^n(A_2)$
birational to $\mathrm{K}^n(A_1)$ is finite up to isomorphism.
Moreover, Theorem \ref{cor:equiv} yields an equivalence
$\Db([A_1/\langle\iota\rangle])\cong\Db([A_2/\langle\iota\rangle])$.
Due to \eqref{eqn:BKR1} and Theorem \ref{cor:equiv},
$\Db(\Km(A))\cong\Db(\Km(B))$ and then $\Km(A)\cong\Km(B)$ (see
\cite{Mu}).\end{proof}

An analogous result for Hilbert schemes of points on K3 surfaces was proved in \cite{Pl}.

\begin{remark}\label{rmk:cont} Observe that, in general, if $A$ and $B$ are abelian
surfaces such that $\Km(A)\cong\Km(B)$, then $\K^n(A)$ and
$\K^n(B)$ are not necessarily birational. Indeed, consider an
abelian surface $A$ such that $A\not\cong\widehat A$ and
$\NS(A)=\langle H\rangle$ with $H^2=6$. Obviously
$\Db(A)\cong\Db(\widehat A)$. Due to Theorem \ref{cor:equiv},
$\Db(\Km(A))\cong\Db(\Km(\widehat A))$ and
$\Km(A)\cong\Km(\widehat A)$ (\cite{Mu}). On the other hand,
Namikawa (\cite[Sect.\ 5]{Na}) proved that $\K^2(A)$ and
$\K^2(\widehat A)$ are not birational.

Furthermore, Example \ref{ex:geom} yields very explicit examples of
isogenous abelian surfaces $A$ and $B$ such that
$\Db(A)\not\cong\Db(B)$. In particular $\Km(A)\not\cong\Km(B)$ and
$\K^n(A)$ is not birational to $\K^n(B)$ for any positive integer
$n$.\end{remark}

\section{Derived categories of twisted Kummer surfaces}\label{sec:brauer}

In this section we prove Theorem \ref{thm:main2} which relates the
existence of equivalences between the twisted derived categories of
two Kummer surfaces and the existence of Hodge isometries between
the generalized transcendental lattices of the corresponding abelian
surfaces. More precisely, in Section \ref{subsec:Br} we introduce
the Brauer group of smooth projective varieties and we prove a
preliminary result (Lemma \ref{lemma:brK3ab}) which will be used in
the proof of Theorem \ref{thm:main2} contained in Section
\ref{subsec:proofsec}. We also discuss a geometric example and an
application to the problem of determining the number of possible
twisted Kummer structures on a twisted K3 surface (respectively in
Sections \ref{ex:geom} and \ref{subsec:num}).

\subsection{Brauer groups and twisted sheaves}\label{subsec:Br} Recall that the {\it Brauer group} $\Br(X)$ of a smooth projective
variety $X$ is the torsion part of $H^2(X,\mathcal{O}_X^*)$ in the
analytic topology (or, equivalently,
$H^2_{\mathrm{\acute{e}t}}(X,\ko_X^*)$ in the \'{e}tale topology)
(see \cite{C1,dJ}).

Assume that $X$ is either a K3 or an abelian surface. It is known
that any $\alpha\in\Br(X)$ is determined (not uniquely) by some
$B\in H^2(X,\QQ)$ (see Chapter 1 of \cite{C1} for the case of K3
surfaces and use a similar argument to deal with abelian surfaces).
This follows from the fact that $H^2(X,\ZZ)$ is unimodular and
$H_1(X,\ZZ)$ is torsion free. More precisely, we deduce the
existence of natural isomorphisms $\Br(X)\cong
T(X)^\vee\otimes\mathbb{Q}/\mathbb{Z}\cong{\rm
Hom}(T(X),\mathbb{Q}/\mathbb{Z})$ and for any $t\in T(X)$,
$\alpha:t\longmapsto t\cdot B\pmod{\ZZ}$, where ``$\cdot$'' is the
cup-product. From this we get a surjective map
\[
\kappa_X:H^2(X,\QQ)\longrightarrow\Br(X).
\]

\begin{lem}\label{lemma:brK3ab}  If $A$ is an abelian surface, there exists an
isomorphism $\Theta_A:\Br(A)\rightarrow\Br(\Km(A))$.\end{lem}

\begin{proof} The K3 surface $\Km(A)$ is the crepant resolution of
$\K(A)=A/\langle\iota\rangle$. Hence there exists a rational map
$\pi:A\dashrightarrow\Km(A)$. Furthermore, as it was observed in
Remark 2 of \cite{Ni1} (see also \cite[Sect.\ 4]{Mo}), the
homomorphism $\pi_{*}$ induces a Hodge isometry
\begin{eqnarray}\label{eqn:trasc}
\pi_{*}:T(A)(2)\longrightarrow T(\Km(A)).
\end{eqnarray}
(Recall that, given a lattice $L$ with quadratic form $b_L$, the
lattice $L(m)$, with $m\in\ZZ$, coincides with $L$ as a group but
its quadratic form $b_{L(m)}$ is such that
$b_{L(m)}(l_1,l_2)=mb_L(l_1,l_2)$, for any $l_1,l_2\in L$.)

In particular, we get a natural morphism
$\Xi:H^2(A,\QQ)\rightarrow T(\Km(A))\otimes\QQ$ defined by
\begin{eqnarray}\label{eqn:Theta}
\Xi:B\longmapsto\frac{\pi_{*}(p(B))}{2},
\end{eqnarray}
where $p:H^2(A,\QQ)\rightarrow T(A)\otimes\QQ$ is the orthogonal
projection. This yields a morphism
$\Theta_A:\Br(A)\longrightarrow\Br(\Km(A))$ of Brauer groups defined
by the commutative diagram
\begin{eqnarray}\label{eqn:diaTheta}
\xymatrix{H^2(A,\QQ)\ar[r]^{\Xi}\ar@{>>}[d]_{\kappa_A}&
T(\Km(A))\otimes\QQ\ar[d]^{\kappa_{\Km(A)}}\\
\Br(A)\ar[r]^{\Theta_A}&\Br(\Km(A)).}
\end{eqnarray}
Observe that $\Theta_A$ is well-defined because, obviously, the
restriction $\kappa_{\Km(A)}|_{T(\Km(A))\otimes\QQ}$ is still
surjective. An easy check then shows that $\Theta_A$ is an
isomorphism.\end{proof}

Any $\alpha\in\Br(X)$ can be represented by a \v{C}ech 2-cocycle on
an analytic cover $\{U_i\}_{i\in I}$ of $X$ using sections
$\alpha_{ijk}\in\Gamma(U_i\cap U_j\cap U_k,\ko^*_X)$. An
\emph{$\alpha$-twisted coherent sheaf} $\kf$ is a pair
$(\{\kf_i\}_{i\in I},\{\varphi_{ij}\}_{i,j\in I})$, where $\kf_i$ is
a coherent sheaf on $U_i$ and $\varphi_{ij}:\kf_j|_{U_i\cap
U_j}\to\kf_i|_{U_i\cap U_j}$ is an isomorphism such that
$\varphi_{ii}=\mathrm{id}$, $\varphi_{ji}=\varphi_{ij}^{-1}$ and
$\varphi_{ij}\circ\varphi_{jk}\circ\varphi_{ki}=\alpha_{ijk}\mathrm{id}$.
Given $\alpha\in\Br(X)$, we denote by $\coh(X,\alpha)$ the abelian
category of $\alpha$-twisted coherent sheaves on $X$ while
$\Db(X,\alpha):=\Db(\coh(X,\alpha))$ is the bounded derived category
of $\coh(X,\alpha)$ (see \cite{C1} and \cite{HS2} for details).

If $X$ and $Y$ are smooth projective varieties and $\alpha\in\Br(X)$
while $\beta\in\Br(Y)$, an equivalence
$\Phi:\Db(X,\alpha)\rightarrow\Db(Y,\beta)$ is a {\it twisted
Fourier-Mukai equivalence} if and only if it there is an isomorphism
as in \eqref{eqn:FMtr} whose kernel $\ke$ is in $\Db(X\times Y,
\alpha^{-1}\boxtimes\beta)$ (see also \cite{CS}).

As in \cite{HS1}, a \emph{twisted variety} is a pair $(X,\alpha)$,
where $X$ is a smooth projective variety and $\alpha\in\Br(X)$. An
isomorphism $f:(X,\alpha)\rightarrow (Y,\beta)$ of the twisted
varieties $(X,\alpha)$ and $(Y,\beta)$ is an isomorphism
$f:X\rightarrow Y$ such that $f^*\beta=\alpha$.

\subsection{Proof of Theorem \ref{thm:main2}}\label{subsec:proofsec} First
of all, observe that, if $X$ is either a K3 or an abelian surface
and $\alpha\in\Br(X)$, the lattice $T(X,\alpha):=\ker(\alpha)$
inherits from $T(X)$ a weight-two Hodge structure. Secondly, if
$\Theta_{A_i}:\Br(A_i)\rightarrow\Br(\Km(A_i))$ is the isomorphism
in Lemma \ref{lemma:brK3ab}, the isometry
$\pi_{i*}:T(A_i)(2)\rightarrow T(\Km(A_i))$ defined in
\eqref{eqn:trasc} yields a Hodge isometry
$$f_i:T(A_i,\alpha)(2)\longrightarrow
T(\Km(A_i),\Theta_{A_i}(\alpha)),$$ for any $\alpha\in\Br(A_i)$
and $i\in\{1,2\}$.

Proposition 4.7 in \cite{Hu}, originally proved for K3 surfaces,
works perfectly in the case of abelian surfaces as well. Therefore
if $X$ is either a K3 or an abelian surface, $\alpha\in\Br(X)$ and
$B\in H^2(X,\QQ)$ is such that $\alpha=\kappa_X(B)$, then there
exists a Hodge isometry
\begin{eqnarray}\label{eqn:Huy}
\begin{array}{rcl}
\exp(B):T(X,\alpha)(k)&\longrightarrow&T(X,B)(k)\\
\gamma&\longmapsto&(\gamma,B\wedge\gamma),
\end{array}
\end{eqnarray}
for any $k\in\{1,2\}$. Given $B_i\in H^2(A_i,\QQ)$, let $\widetilde
B_i\in H^2(\Km(A_i),\QQ)$ be such that
$$\Theta_{A_i}(\kappa_{A_i}(B_i))=\kappa_{\Km(A_i)}(\widetilde B_i).$$ Define
$\widetilde\alpha_i:=\Theta_{A_i}(\kappa_{A_i}(B_i))$. If $g:
T(A_1,B_1)\rightarrow T(A_2,B_2)$ is a Hodge isometry, the diagram
\[
\xymatrix{T(A_1,B_1)(2)\ar[rrr]^{g}\ar[ddd]\ar[dr]_{\exp(-B_1)}&&&
T(A_2,B_2)(2)\ar[ddd]\ar[dl]^{\exp(-B_2)}\\
& T(A_1,\kappa_{A_1}(B_1))(2)\ar[d]_{f_1}\ar[r]& T(A_2,\kappa_{A_2}(B_2))(2)\ar[d]^{f_2}&\\
& T(\Km(A_1),\widetilde\alpha_1)\ar[r]& T(\Km(A_2),\widetilde\alpha_2)&\\
T(\Km(A_1),\widetilde B_1)\ar[rrr]^{f}\ar[ur]^{\exp(\widetilde
B_1)}&&& T(\Km(A_2),\widetilde B_2)\ar[ul]_{\exp(\widetilde B_2)}}
\]
commutes and yields a Hodge isometry $f: T(\Km(A_1),\widetilde
B_1)\rightarrow T(\Km(A_2),\widetilde B_2)$. Conversely, since
$\Theta_i$ is an isomorphism (Lemma \ref{lemma:brK3ab}), the same
diagram and remarks show that any Hodge isometry between the
generalized transcendental lattices of $\Km(A_1)$ and $\Km(A_2)$
determined by some $\widetilde B_i\in H^2(\Km(A_i),\QQ)$ induces a
Hodge isometry of the generalized transcendental lattices of $A_1$
and $A_2$ corresponding to $B_i\in H^2(A_i,\QQ)$ such that
$$\kappa_{A_i}(B_i)=\Theta_{A_i}^{-1}(\kappa_{\Km(A_i)}(\widetilde
B_i))\in\Br(A_i).$$

Since the Picard number of $\Km(A_i)$ is greater than 11, the
equivalence between item (i) and item (ii) of Theorem
\ref{thm:main2} follows from \cite[Thm.\ 0.4]{HS1}. Indeed such a
result proves that, for any $B_i\in H^2(\Km(A_i),\QQ)$, there
exists a twisted Fourier-Mukai equivalence
\[
\Db(\Km(A_1),\kappa_{\Km(A_1)}(B_1))\cong\Db(\Km(A_2),\kappa_{\Km(A_2)}(B_2))
\]
if and only if there exists a Hodge isometry
$T(\Km(A_1),B_1)\cong T(\Km(A_2),B_2)$.

Due to what we have just proved, any twisted Fourier-Mukai
equivalence $\Db(\Km(A_1),\alpha_1)\cong\Db(\Km(A_2),\alpha_2)$
induces a Hodge isometry $T(\Km(A_1))\otimes\QQ\cong
T(\Km(A_2))\otimes\QQ$ which yields a Hodge isometry
$T(A_1)\otimes\QQ\cong T(A_2)\otimes\QQ$. Consider the Kuga-Satake
varieties $\KS(A_1)$ and $\KS(A_2)$ associated to the weight-two
Hodge structures on $T(A_1)\otimes\QQ$ and $T(A_2)\otimes\QQ$ (see
Section 4 in \cite{Mo1} for the definition). Theorem 4.3 and Lemma
4.4 in \cite{Mo1} show that, for any $i\in\{1,2\}$,
$$\underbrace{\KS(A_i)\times\cdots\times\KS(A_i)}_{2^{\rho(A_i)}\mbox{ times}}\sim\underbrace{A_i\times\cdots\times A_i}_{8\mbox{ times}},$$ where ``$\sim$'' denotes
an isogeny of abelian varieties. By construction
$\KS(A_1)\sim\KS(A_2)$ and then $A_1^8\sim A_2^8$. In particular,
$A_1$ and $A_2$ are isogenous and this concludes the proof of
Theorem \ref{thm:main2}.

To shorten the notation and according to \cite{HS1}, we introduce
two equivalence relations:

\begin{definition}\label{def:DTequiv} Let $(X_1,\alpha_1)$ and $(X_2,\alpha_2)$ be twisted K3
or abelian surfaces.
\begin{itemize}
\item[(i)] They are {\it $D$-equivalent} if there exists a twisted
Fourier-Mukai equivalence
$$\Phi:\Db(X_1,\alpha_1)\rightarrow\Db(X_2,\alpha_2).$$

\item[(ii)] They are {\it $T$-equivalent} if there exist $B_i\in
H^2(X_i,\QQ)$ such that $\alpha_i=\kappa_{A_i}(B_i)$ and a Hodge
isometry
$$\varphi:T(X_1,B_1)\rightarrow
T(X_2,B_2).$$\end{itemize}\end{definition}

We can now prove the following easy corollary of Theorem
\ref{thm:main2}:

\begin{cor}\label{cor:cons} {\rm (i)} $(\Km(A_1),1)$ is
D-equivalent to $(\Km(A_2),1)$ if and only if $(A_1,1)$ and
$(A_2,1)$ are T-equivalent.

{\rm (ii)} If $(A_1,\alpha_1)$ and $(A_2,\alpha_2)$ are D-equivalent
twisted abelian surfaces, then $(\Km(A_1),\Theta_{A_1}(\alpha_1))$
and $(\Km(A_2),\Theta_{A_2}(\alpha_2))$ are D-equivalent.\end{cor}

\begin{proof} Due to the isomorphism in Lemma \ref{lemma:brK3ab}, (i) follows trivially from Theorem \ref{thm:main2}. The machinery in \cite{HS1} applied to the case of abelian surfaces shows that if $(A_1,\alpha_1)$ and $(A_2,\alpha_2)$ are D-equivalent, then they are T-equivalent as well. Then use Theorem \ref{thm:main2}.\end{proof}

Notice that part (i) of Corollary \ref{cor:cons} is exactly the
analogue of (B) in the introduction.

\begin{remark}\label{rmk:commenti} (i) Due to \cite[Prop.\ 8.1]{HS1},
if $\alpha_j\in\Br(\Km(A_j))$ is non-trivial for any
$j\in\{1,2\}$, then the existence of an equivalence
$\Db(\Km(A_1),\alpha_1)\cong\Db(\Km(A_2),\alpha_2)$ does not imply
that $\Km(A_1)\cong\Km(A_2)$ (see also Example \ref{ex:geom}). This is one of the main differences
with the untwisted case treated by Hosono, Lian, Oguiso and Yau in
\cite{HLOY} (see (A) and (C) in the introduction).

(ii) As suggested by Corollary \ref{cor:cons}, we would expect (ii)
in Theorem \ref{thm:main2} to be equivalent to the existence of a
twisted Fourier-Mukai equivalence
$\Db(A_1,\beta_1)\cong\Db(A_2,\beta_2)$, where $\beta_i\in\Br(A_i)$.
This would lead to a twisted version of (A). Actually this is not
the case. Indeed, since the period map is surjective for abelian
surfaces (\cite{Sh}), one can produce a counterexample to this
expectation by adapting Example 4.11 in \cite{HS1}.

(iii) Let $A_1$ and $A_2$ be two abelian surfaces with
$\NS(A_1)=\langle H_1\rangle$ and $\NS(A_2)=\langle H_2\rangle$. If
there exist $\alpha_1\in\Br(\Km(A_1))$ and
$\alpha_2\in\Br(\Km(A_2))$ such that
$\Db(\Km(A_1),\alpha_1)\cong\Db(\Km(A_2),\alpha_2)$ then
$H_1^2/H_2^2$ is a square in $\QQ$. Indeed, by Theorem
\ref{thm:main2} (and by \cite[Sect.\ 7]{HS1}), if
$\Db(\Km(A_1),\alpha_1)\cong\Db(\Km(A_2),\alpha_2)$ then there
exists an isogeny $\varphi:A_1\rightarrow A_2$ inducing a Hodge
isometry $\varphi^*:H^2(A_2,\QQ)\rightarrow H^2(A_1,\QQ)$ such that
$\varphi^*(H_2)=qH_1$, for some $q\in\QQ$. In particular
$H_2^2=q^2H_1^2$.\end{remark}

\subsection{An explicit example}\label{ex:geom} In this example, we use Theorem
\ref{thm:main2} to establish a connection between the twisted
derived categories of some nice Kummer surfaces with Picard number
2. Recall that the lattices $U$ and $U(n)$ are the free abelian
group $\ZZ\oplus\ZZ$ endowed respectively with the quadratic forms
represented by the matrices
\[
\left(\begin{array}{cc} 0 & 1\\ 1 &
0\end{array}\right)\;\;\;\mbox{and}\;\;\;\left(\begin{array}{cc} 0 & n\\
n & 0\end{array}\right).
\]

Let $A$ be an abelian surface such that $\NS(A)\cong U(n)$, for
some positive integer $n$. We first show that there exist two
elliptic curves $E$ and $F$ and a subgroup $C_n\cong\ZZ/n\ZZ$ of
$E\times F$ such that either $A\cong(E\times F)/C_n$ or
$\widehat{A}\cong(E\times F)/C_n$.

To see this, let us first observe that, since $\NS(A)\cong U(n)$,
the transcendental lattice $T(A)$ is isometric to $U(n)\oplus U$.
Indeed for any abelian surface $A$, $H^2(A,\ZZ)$, endowed with the
cup-product, is isometric to the lattice $U\oplus U\oplus U$ (see
\cite{Mo} for more details).

We choose a basis $\langle e_1,e_2,f_1,f_2\rangle=U\oplus
U(n)\hookrightarrow U^3$, an isometry $\varphi:H^2(A,\ZZ)\rightarrow
U^3$ and $c\in\CC$ such that
\begin{eqnarray}\label{eqn:period}
\varphi(c\sigma_A)=e_1-n\omega_1\omega_2e_2+\omega_1f_1+\omega_2f_2,
\end{eqnarray}
where $H^{2,0}(A)=\langle\sigma_A\rangle$. We define in $\CC$ the
lattices $\Gamma_1=\ZZ+\omega_1\ZZ$ and $\Gamma_2=\ZZ+\omega_2\ZZ$
and the elliptic curves $E:=\CC/\Gamma_1$ and $F:=\CC/\Gamma_2$.
Notice that, since $T(A)\cong U(n)\oplus U$ and $\sigma_A^2=0$,
the numbers $1$, $\omega_1$, $\omega_2$ and $\omega_1\omega_2$ are
linearly independent over $\QQ$. So, in particular, $E$ and $F$
are not isogenous. If $H_1(E\times
F,\ZZ)=\langle\gamma_1,\gamma_2,\delta_1,\delta_2\rangle$, then we
consider the subgroup $C_n$ of $E\times F$ such that
\[
H_1((E\times
F)/C_n,\ZZ)=\left\langle\frac{\gamma_1+\delta_1}{n},\gamma_2,\delta_1,\delta_2\right\rangle.
\]
Let $S:=(E\times F)/C_n$. In terms of the dual bases of the bases of
$H_1(E\times F,\ZZ)$ and  $H_1(S,\ZZ)$ just described, we write
$H^1(S,\ZZ)=\langle\de z_1,\de z_2,\de w_1,\de w_2\rangle$ and
$H^1(E\times F,\ZZ)=\langle\de x_1,\de x_2,\de y_1,\de y_2\rangle$.
If $\pi:E\times F\rightarrow S$ is the natural surjection, the map
$\theta:=\pi^*:H^1(S,\ZZ)\rightarrow H^1(E\times F,\ZZ)$ is such
that:
\begin{eqnarray}\label{eqn:theta}
\begin{array}{l}
\theta(\de z_1)=n\de x_1,\\
\theta(\de z_2)=\de x_2,\end{array}\;\;\;\;\;\;\;\;\;\;\begin{array}{l}\theta(\de w_1)=-\de x_1+\de y_1,\\
\theta(\de w_2)=\de y_2.
\end{array}
\end{eqnarray}
Observe that $\NS(E\times F)=\langle\de x_1\wedge\de x_2,\de
y_1\wedge\de y_2\rangle$. Furthermore, due to the properties in
\eqref{eqn:theta} which characterize the morphism
$\stackrel{2}{\wedge}\theta:H^2(S,\ZZ)\rightarrow H^2(E\times
F,\ZZ)$ and due to the fact that $\stackrel{2}{\wedge}\theta$
preserves the Hodge structures on $H^2(S,\ZZ)$ and $H^2(E\times
F,\ZZ)$,
{\setlength\arraycolsep{2pt}
\begin{eqnarray*}
\NS(S)&=&\langle\de z_1\wedge\de z_2,n\de w_1\wedge\de w_2+\de
z_1\wedge\de w_2\rangle,\\
T(S)&=&\langle\de z_1\wedge\de w_1,\de z_2\wedge\de w_2,\de
z_1\wedge\de w_2,-n\de w_1\wedge\de z_2+\de z_1\wedge\de w_2\rangle.
\end{eqnarray*}}In particular, $\NS(S)\cong U(n)$ and $T(S)\cong U\oplus U(n)$.

Consider the two cohomology classes
{\setlength\arraycolsep{2pt}
\begin{eqnarray*}
\sigma_{E\times F}&=&\de x_1\wedge\de y_1+\omega_2\de x_1\wedge\de
y_2+\omega_1\de x_2\wedge\de y_1+\omega_1\omega_2\de x_2\wedge\de
y_2;\\
\sigma_S&=&\de z_1\wedge\de w_1+\omega_2\de z_1\wedge\de
w_2+\omega_1(n\de w_1\wedge\de z_2-\de z_1\wedge\de
w_2)+n\omega_1\omega_2\de z_2\wedge\de w_2.
\end{eqnarray*}}Obviously, $\sigma_{E\times F}\in T(E\times F)\otimes\CC$ and
$\sigma_S\in T(S)\otimes\CC$. Since $\langle\sigma_{E\times
F}\rangle=H^{2,0}(E\times F)$ and since an easy calculation shows
that $\stackrel{2}{\wedge}\theta(\sigma_S)=n\sigma_{E\times F}$,
$\langle\sigma_S\rangle=H^{2,0}(S)$. This implies that, due to
\eqref{eqn:period}, there exists an isometry
$\eta:H^2(S,\ZZ)\rightarrow U^3$ such that
$\eta^{-1}\circ\varphi:H^2(A,\ZZ)\rightarrow H^2(S,\ZZ)$ is a Hodge
isometry (see \cite{Sh}). The Torelli Theorem for abelian surfaces
shows that either $A\cong(E\times F)/C_n$ or
$\widehat{A}\cong(E\times F)/C_n$.

Observe that, since $\NS(A)\cong U(n)$, the abelian surface $A$ is
principally polarized if and only if $n=1$. This means that, if
$n\neq 1$, $(E\times F)/C_n$ and its dual are not isomorphic.
Furthermore, $A$ and $E\times F$ are isogenous but $T(A)\not\cong
T(E\times F)$. Therefore, due to Theorem \ref{thm:orlov1},
$\Db(A)\not\cong\Db(E\times F)$. This proves that there are
isogenous abelian surfaces whose derived categories of (untwisted)
coherent sheaves are not equivalent (see Remark \ref{rmk:cont}).

Choose the standard basis $\{g_1,g_2,k_1,k_2\}$ for $U\oplus U$. Due
to the explicit description of $T(A)$ that we have previously given,
it is straightforward to see that there exists an inclusion
$i_1:T(A)\rightarrow U\oplus U$ where $i_1(e_j)=g_j$
($j\in\{1,2\}$), $i_1(f_1)=nk_1$ and $i_1(f_2)=nk_2$. Let
$\sigma:=i(\sigma_A)\in U^2\otimes\CC$. Due to \eqref{eqn:period},
we can write
$\sigma=g_1-n\omega_1\omega_2g_2+n\omega_1h_1+\omega_2h_2$.

Consider in $\CC$ the lattice $\Gamma_3=\ZZ+n\omega_1\ZZ$ and the
elliptic curve $E_1:=\CC/\Gamma_3$. Of course, $E$ and $E_1$ are
isogenous. Reasoning as before and using the surjectivity of the
period map and the Torelli Theorem for abelian surfaces (\cite{Sh}),
we get an isometry $\varphi_1:T(E_1\times F)\rightarrow U^2$ fitting
in the following commutative diagram:
\[
\xymatrix{0\ar[dr] & & T(E_1\times
F)\ar[dd]^{\varphi_1}\ar[rd]^{\alpha} & & 0 & \\  &
T(A)\ar[dr]^{i_1}\ar[ru]^{i} &  & \ZZ/n\ZZ\ar[rd]\ar[ru] &  \\
0\ar[ru] &
 & U\oplus U\ar[ru]^{\alpha_1} &  & 0.}
\]
Of course, $i$ preserves the Hodge structures and
$\alpha\in\Br(E_1\times F)$. Proposition 4.7 in \cite{Hu} yields
$B\in H^2(E\times F,\QQ)$ such that $(A,1)$ and $(E\times
F,\kappa_{E\times F}(B))$ are T-equivalent. By Theorem
\ref{thm:main2}, there exist $\beta\in\Br(\Km(E_1\times F))$ of
order $n$ and a twisted Fourier-Mukai equivalence
$\Db(\Km(A))\cong\Db(\Km(E_1\times F),\beta)$.

\subsection{The number of twisted Kummer structures}\label{subsec:num} As an easy corollary of Lemma \ref{lemma:brK3ab}, we get a
surjective map
\[
\Psi:\{\mbox{Twisted abelian
surfaces}\}/\mathrm{isom}\longrightarrow\{\mbox{Twisted Kummer
surfaces}\}/\mathrm{isom}
\]
which sends the isomorphism class $[(A,\alpha)]$ to the isomorphism
class $[(\Km(A),\Theta_A(\alpha))]$. The main result in \cite{HLOY}
proves that the preimage of $[(\Km(A),1)]$ is finite, for any
abelian surface $A$ and $1\in\Br(A)$ the trivial class (see
\cite[Thm.\ 0.1]{HLOY}). On the other hand \cite{HLOY} shows that
the cardinality of the preimages of $\Psi$ can be arbitrarily large.
This answered an old question by Shioda. Namely, there can be many
non-isomorphic (untwisted) abelian surfaces giving rise to
isomorphic (untwisted) Kummer surfaces (a partial result in this
direction is also contained in \cite{GH}). This is usually rephrased
saying that on a K3 surface one can put many non-isomorphic
(untwisted) Kummer structures.

Using Theorem \ref{thm:main2}, the picture in \cite{HLOY} can be
completely generalized to the twisted case.

\begin{prop}\label{prop:twnum} {\rm (i)} For any twisted Kummer surface
$(\Km(A),\alpha)$, the preimage $\Psi^{-1}([(\Km(A),\alpha)])$ is
finite.

{\rm (ii)} For positive integers $N$ and $n$, there exists a twisted
Kummer surface $(\Km(A),\alpha)$ with $\alpha$ of order $n$ in
$\Br(\Km(A))$ and such that $|\Psi^{-1}([(\Km(A),\alpha)])|\geq
N$.\end{prop}

\begin{proof} Suppose that
$\Psi([(A_1,\alpha_1)])=\Psi([(A_2,\alpha_2)])=[(\Km(A),\alpha)]$,
i.e.\ suppose that there exists an isomorphism
$f:\Km(A)\isomor\Km(A_i)$ such that
$f^*\Theta_{A_i}(\alpha_i)=\alpha$. In particular,
$$\Db(\Km(A_1),\Theta_{A_1}(\alpha_1))\cong\Db(\Km(A_2),\Theta_{A_2}(\alpha_2)).$$

Due to Theorem \ref{thm:main2}, the proof of (i) amounts to show
that, up to isomorphisms, there are finitely many T-equivalent
twisted abelian surfaces $(A',\beta)$ such that
$\Psi([(A',\beta)])=[(\Km(A),\alpha)]$. Since, up to isomorphisms,
there are just finitely many abelian surfaces $A'$ with
$\Db(A')\cong\Db(A)$ (\cite[Prop.\ 5.3]{BM}), we can just fix $A'$
with such a property and show that, up to isomorphisms, there exists
a finite number of $\beta'\in\Br(A')$ such that $(A',\beta)$ and
$(A',\beta')$ are T-equivalent. But this is the content of
\cite[Prop.\ 3.4]{HS1} for the case of abelian surfaces.

Applying the results in \cite{Og} and \cite{St} to abelian surfaces,
we see that, for any positive integer $N$, there exist $N$
non-isomorphic abelian surfaces $A_1,\ldots,A_N$ such that
$\Db(A_i)\cong\Db(A_j)$ ($i,j\in\{1,\ldots,N\}$). Due to Theorem
\ref{thm:orlov1}, for any $i\in\{2,\ldots,N\}$, there is a Hodge
isometry $$g_i:T(A_1)\rightarrow T(A_i).$$ Take $B_1\in
T(A_1)\otimes\QQ$ such that $\alpha_1:=\kappa_{A_1}(B_1)$ and
$\Theta_{A_1}(\alpha_1)$ are not trivial in $\Br(A_1)$ and
$\Br(\Km(A_1))$ respectively. We can also choose $\alpha_1$ such
that the order of $\Theta_{A_1}(\alpha_1)$ is $n$ in
$\Br(\Km(A_1))$. Then, for any $i\in\{2,\ldots,N\}$, define
$\alpha_i:=\kappa_{A_i}(g_i(B_1))$. Obviously, $(A_i,\alpha_i)$ and
$(A_j,\alpha_j)$ are T-equivalent when $i,j\in\{1,\ldots,N\}$.
Theorem \ref{thm:main2} immediately implies that
$(\Km(A_i),\Theta_{A_i}(\alpha_i))$ and
$(\Km(A_j),\Theta_{A_j}(\alpha_j))$ are D-equivalent.

For any $i\in\{2,\ldots,N\}$, the isometry $g_i$ induces a Hodge
isometry $f_i:T(\Km(A_1))\rightarrow T(\Km(A_i))$ which (due to
\cite[Thm.\ 1.14.4]{Ni}) extends to a Hodge isometry
$h_i:H^2(\Km(A_1),\ZZ)\rightarrow H^2(\Km(A_i),\ZZ)$. The Torelli
Theorem yields an isomorphism $\varphi_i:\Km(A_1)\rightarrow
\Km(A_i)$ such that
$\varphi_i^*(\Theta_{A_i}(\alpha_i))=\Theta_{A_1}(\alpha_1)$
(possibly changing $\alpha_i$ with $\alpha_i^{-1}$), for any
$i\in\{2,\ldots,N\}$. This concludes the proof of (ii).\end{proof}

In other words, Proposition \ref{prop:twnum} shows that on a twisted
K3 surface we can put just a finite number of non-isomorphic
\emph{twisted Kummer structures}. Nevertheless, such a number can be
arbitrarily large even when the twist is non-trivial and has any
possible order.

\medskip

{\small\noindent{\bf Acknowledgements.}  The final version of this
paper was written during the author stay at the Max-Planck-Institut
f\"{u}r Mathematik (Bonn). The results in this paper benefit of many
interesting discussions with Daniel Huybrechts and Bert van Geemen.
It is a great pleasure to thank them for their help and
suggestions.}

\end{document}